\documentclass[review]{elsarticle}

\usepackage{hyperref}

\usepackage{float}
\usepackage{placeins}
\usepackage{mathtools}
\usepackage{amsfonts}
\usepackage{amsmath}
\usepackage{amssymb}
\usepackage{graphicx}
\usepackage{float}
\usepackage{verbatim}
\usepackage{array}
\usepackage{multirow}
\usepackage{dcolumn}
\usepackage{color}

\usepackage{setspace}
\usepackage{diagbox}
\usepackage{adjustbox} 
\usepackage{makecell}
\usepackage{mathtools}
\usepackage{comment}

\usepackage{geometry}
\journal{Journal of \LaTeX\ Templates}

\begin{document}

\begin{frontmatter}


\title{Solving Multistage Stochastic Linear Programming via Regularized Linear Decision Rules: An Application to Hydrothermal Dispatch Planning}


\author{Felipe Nazare and Alexandre Street*}
\address{LAMPS PUC-Rio,\\ Electrical Engineering Department}
\fntext[myfootnote]{Felipe Nazare and Alexandre Street are with the electrical engineering Department at the Pontifical Catholic University of Rio de Janeiro (PUC-Rio). Alexandre Street is also the co-founder and research director at the Laboratory of Applied Mathematical Programming and Statistics (LAMPS) PUC-Rio.}

\author[mymainaddress]{Pontifical Catholic University of Rio de Janeiro}

\cortext[mycorrespondingauthor]{Corresponding author: Alexandre Street}
\ead{street@puc-rio.br}


\begin{abstract}
The solution of multistage stochastic linear problems (MSLP) represents a challenge for many application areas. Long-term hydrothermal dispatch planning (LHDP) materializes this challenge in a real-world problem that affects electricity markets, economies, and natural resources worldwide. No closed-form solutions are available for MSLP and the definition of non-anticipative policies with high-quality out-of-sample performance is crucial. Linear decision rules (LDR) provide an interesting simulation-based framework for finding high-quality policies for MSLP through two-stage stochastic models. In practical applications, however, the number of parameters to be estimated when using an LDR may be close to or higher than the number of scenarios of the sample average approximation problem, thereby generating an in-sample overfit and poor performances in out-of-sample simulations. In this paper, we propose a novel regularized LDR to solve MSLP based on the AdaLASSO (adaptive least absolute shrinkage and selection operator). The goal is to use the parsimony principle as largely studied in high-dimensional linear regression models to obtain better out-of-sample performance for LDR applied to MSLP. Computational experiments show that the overfit threat is non-negligible when using classical non-regularized LDR to solve the LHDP, one of the most studied MSLP with relevant applications. Our analysis highlights the following benefits of the proposed framework in comparison to the non-regularized benchmark: 1) significant reductions in the number of non-zero coefficients (model parsimony), 2) substantial cost reductions in out-of-sample evaluations, and 3) improved spot-price profiles.


\end{abstract}


\end{frontmatter}







\section{Introduction}

Whenever a decision process of one period affects the optimal decisions of the subsequent ones, it gives rise to a dynamic problem, which in the optimization field is known as a multistage problem. When this process is affected by uncertainty that is gradually observed and decisions can also be dynamically adapted to the uncertainty process through a linear programming problem, a multistage stochastic linear problem (MSLP) can be defined.

The solution of a MSLP relies on the estimation of a decision policy, i.e., a rule defining the vector of decision variables to be implemented throughout the stages (periods) as a function of the history of the uncertain parameters representing the state of the system. For the general MSLP, there is no closed form available. Furthermore, the combination of a high number of stages and high dimension of state and uncertainty vectors results in an enormous range of possible paths, preventing the use of scenario-tree approaches (see \cite{birge2011}). Thus, the predominant literature exploring MSLP as a simulation-based tool for practical decision-making under uncertainty widely focuses on sample average approximation (SAA) techniques (see \cite{saa99} and \cite{SPnewsletter2}). This phenomenon has been named in the literature as \textit{the curse of dimensionality}. If we consider a stage-wise independent uncertainty affecting the right-hand side of a linear problem, the so-called stochastic dual dynamic programming (SDDP) method provides us with the optimal policy for a given sample of scenarios. This concept was proposed in \cite{Pereira1991} applied to long-term hydrothermal dispatch planning (LHDP), and its convergence was proved in \cite{Philpott2008}. Just to mention a few, more recently, SDDP was applied in \cite{Street2017} to consider security criteria in the dispatch problem through a hybridization with robust optimization and in \cite{Tony2018} to address the short-term dispatch with storage in a system affected by wind uncertainty. 

Alternatively, a decision policy for a MSLP can also be defined through a function of observed data, say, ${\xi}_{[1:t]}$, where $t$ refers to the stage, affecting the problem until the moment each vector of decisions, ${x}_t$, must be defined, i.e., ${x}_t = {\Psi}(\xi_{[1:t]})$. Approximation approaches consider the subset of policies defined as a linear combination of functions applied to the uncertainty data. This subset of decision rules is also known as linear decision rules (LDRs), oftentimes referred to as affine rules or affine policies in cases where an affine functional space is used. We refer to \cite{Bodur_2018} and references therein for a recent work. This method first appeared in \cite{Holt1955}, applying it in an employment workforce and factory production. However, due to its flexibility in solving different stochastic problems, the approach was easily accommodated to other linear problems, being widely used in reservoir management in a large set of publications such as  \cite{Revelle1969,Revelle1970,Resources1970,Eastman1973,Gundelach1975,Revelle1975}.

The subject of LDR in decision-making under uncertainty is still under development, as observed in \cite{Gauvin2017}. In this {work}, the authors aim to minimize worst-case scenarios of reservoir limit violations through a robust approach, approximating the total outflow curve to affine functions of the uncertainty data. Besides reservoir management problems, the simple adaptation of this method allowed its application in many other studies, varying from the solution of financial assets allocation problem \cite{Calafiore_2008} to the management of heat and power generation \cite{Zugno2014}, power system operation \cite{Braaten2016}, power generation capacity planning \cite{Dominguez2016} among other relevant works. In terms of hydropower scheduling, the application of LDR in \cite{Egging2017} aims the definition of a feasible hydropower dispatch for a price and inflow uncertainty set, allowing for a reduction in the computational effort to solve the problem while obtaining gains in expected profits. It is relevant to mention that LDR has been playing a relevant role in addressing two-stage and multi-stage robust models in the recent literature (see \cite{Lorca2017} and references therein). This is because the LDR allows for a simplified representation of recourse actions, avoiding the complexity of multi-level optimization that robust models generally impose.

So, while SDDP-based approaches limit the uncertainty model, the type of variables that can be considered, and require convexity, the use of LDR limits the format of the functional used to define the decision rule. However, it is worth mentioning that the framework of LDR is not limited to linear or affine functions of the uncertainties. Rather, any convenient linear functional space of interest mapping the set of uncertainty vectors onto the set of states can be considered and the {problem} linearity would still be preserved\footnote{{Note that defining a LDR within a given linear space of nonlinear functions (trigonometric, polynomial, etc.) is equivalent to the problem of finding the best coefficients of a basis for this space, which is a linear problem.}}. Additionally, as shown in \cite{Bodur_2018}, the adoption of a LDR approach allows writing the MSLP as a two-stage stochastic model (with multiple periods). In this case, we would be inheriting all properties and convergence results (e.g., SAA results \cite{saa99,shapiro2021lectures}), and the methods and algorithms (such as Benders decomposition \cite{sikora2021benders} and Progressive Hedging \cite{nikzad2021matheuristic}) for this widely studied class of models. Regarding its flexibility, this class of models allows us to address MSLP based on uncertainty processes with any kind of non-linear time-dependency structure and non-Gaussian distributions (such as those typically affecting power systems -- see \cite{GAS} and \cite{bodin2020scoredrivenmodels}) through completely exogenously-generated scenarios, to consider the broader class of mixed-integer linear problems, and to consider data-driven methods (see \cite{velloso2019two}). In this context, it is worth emphasizing that both approaches have pros and cons and deserve further studies, as recently demonstrated in the comparisons provided in \cite{Bodur_2018}. 

\emph{In this work, we work concentrate on providing a new {empirical and practical contribution to the application of LDRs to solve MSLPs}.} Within this context, the definition of a policy requires the estimation of many coefficients under a limited amount of data. This is especially critical in the LHDP problem, where a huge number of coefficients would be needed to account for each stage, individual reservoirs, different components of a functional basis, and lags of observed data in the decision rule. Hence, in practical applications, we oftentimes need to handle more coefficients than the number of scenarios to keep the model tractable. In this situation, we should expect from LDRs a certain degree of in-sample overfit and, consequently, poor out-of-sample performances, i.e., under unseen scenarios. This is a well-known fact from high-dimensional statistics and machine learning problems (we refer to \cite{zou2006adaptive,giraud2014,ciuperca_adaptive_2016,blanchet2019robust} and \cite{medeiros2020} as a relevant and updated literature on the subject). 

Interestingly, when applied to a MSLP, the two-stage LDR identification problem resembles the estimation of a linear regression model, $Y=X'\theta+\varepsilon$, based on a more general loss function, $l(X,Y;\theta)$, related to the specific objective function of the problem. In this context, the overfitting issue and the related poor out-of-sample performance of non-parsimonious {regression} models (with the number of coefficients approaching or exceeding the number of in-sample observations) have been addressed by the use of regularization procedures (\cite{medeiros2020}). In particular, the \emph{adaptive least absolute shrinkage and selection operator} (AdaLASSO) has been playing a key role in the statistical and machine learning literature \cite{zou2006adaptive, ciuperca_adaptive_2016,  Tibshirani1996, medeiros2020}. The interested reader is referred to \cite{zou2006adaptive} and \cite{giraud2014} for further details. More recently, in \cite{blanchet2019robust}, the capacity to hedge against unseen data (out-of-sample) of regularized regression models, {such as those based on the LASSO}, was connected to data-driven distributionally robust optimization problems with Wasserstein metric. {In \cite{bertsimas2019book} (see chapter 2), linear regression with LASSO and other regularization metrics were also connected to robust regression models.} { Although such theoretical results cannot be generalized for a general loss function, previously} reported (empirical and theoretical) results regarding the improved out-of-sample performance of regularized { regression} models provide us with strong { and insightful} evidence that MSLP addressed via LDR would also benefit from regularization schemes.

\subsection{Hypothesis, objectives, and contributions}
The main hypothesis tested in this work is the following: \emph{non-regularized LDR overfits in-sample data and more parsimonious LDR-based policies estimated with AdaLASSO-regularization can provide relevant cost savings under unseen out-of-sample data}. {Thus, the objectives of this paper are: 1) to test the aforementioned hypothesis, 2) to propose a novel regularized LDR for MSLP based on the AdaLASSO, and 3) to study the regularization effects using the LHDP problem. The goal of the second objective is to improve the out-of-sample performance of a LDR estimated with a finite sample of size $N$.} 

To achieve the proposed objectives, this paper provides the following contributions to the state-of-the-art literature: 

\begin{itemize}
\setlength\itemsep{1em}
 \item {We raise awareness of the in-sample overfit threat to which LDRs are subject when applied to solve MSLP. We showcase significant over costs due to overfit and characterize other detrimental effects due to non-regularized policies on both primal and dual (spot prices) solutions of the LHDP.}
 
 
  \item {We show strong evidence that the proposed regularization process is an interesting avenue to robustify the estimation process of LDR-based policies against overfit. In this sense, as our approach provides an out-of-sample cost lower or equal to the out-of-sample cost of the non-regularized LDR, we also improve the SAA method (\cite{saa99}) applied to solve MSLP based on LDR by offering a lower or equal upper bound and an associated improved implementable solution. We exemplify the benefits of the proposed AdaLASSO-regularized LDR by studying one of the most relevant MSLP found in the related literature, the LHDP, using a polynomial basis of functions.}

 
\end{itemize}
 
To quantify the benefits of our proposed regularization scheme, we use a two-step estimation-simulation procedure. In the first step, the LDR coefficients are estimated by a two-stage linear programming problem, in which an AdaLASSO-regularization term imposes the trade-off between parsimony and cost reduction for the in-sample scenarios. In the second step, the state-target tracking problem proposed in \cite{Bodur_2018} is used to find implementable decisions that follow as much as possible the estimated LDR for a large set of out-of-sample scenarios. Computational experiments showcase that the parsimony induced by the regularization procedure significantly reduces the number of coefficients different from zero, shrinks the remaining nonzero coefficients, and provides {significant out-of-sample economic and operational gains} for the identified policy solving the LHDP problem. 

{Finally, it is worth mentioning that our simple yet novel idea of regularized LDRs provides us with an effective machine-learning-based approach to improve the out-of-sample performance of two-stage LDRs proposed in \cite{Bodur_2018} that does not increase the complexity of the original problem. Indeed, according to \cite{Kuhn2018} and \cite{bertsimas2022two}, the general Wasserstein-based data-driven distributionally robust counterpart of a polynomial-time solvable linear optimization problem with RHS uncertainty is an NP-hard problem. And based on the results of \cite{gamboa2021}, it is actually a ``practically hard NP-hard" problem\footnote{We call a practically hard NP-hard problem a problem that is not only hard for the worst-case instance, but for practical medium-sized instances of interest in realistic applications.}. Therefore, our proposed approach can also be seen as a tractable heuristic to the task of robustification against unseen data that do not add an NP-hard bilevel layer on top of the already difficult original LHDP problem.}

The rest of this paper is structured as follows: Section \ref{Hydrothermal Dispatch Problem} presents the LDR-based policy estimation and evaluation procedures applied to the long-term hydrothermal dispatch planning problem. Section \ref{Regularized_AffinePol} describes the new estimation model with AdaLASSO regularization and the calibration scheme for the regularization parameter. Section \ref{Quantitative_Analysis} depicts the case study data and the main results. Finally, Section \ref{Conclusion} presents the main conclusions of this work and future research avenues.

\section{Hydrothermal Dispatch Problem with Linear Decision Rules}
\label{Hydrothermal Dispatch Problem}

Uncertainty is an inherent feature of any planning activity. The impact of such uncertainties in the decision processes requires well-structured and robust planning processes. This is especially evident in the operation of hydrothermal power systems, as poor management of the water resources may impose long and severe economic and social losses. { The LHDP problem is one of the most studied MSLP. Its representative mathematical structure (MSLP with long horizons and high-dimensional state space and uncertainty) and real-world economic and social relevance (management of one of the most relevant natural resources used to generate energy) motivated a myriad of academic papers (see the literature review in \cite{street2020assessing}) and industry applications in the last three decades.\footnote{{Interestingly, the development of the stochastic programming area was significantly affected by the relevance and challenge of this problem. We refer the interested reader to \cite{SPnewsletter2}.}}} 

In the LHDP problem, independent system operators (ISOs) must manage { the water stored in hydro reservoirs through time to dispatch the best mix of generators to supply demand at every stage while minimizing the long-run expected operating cost. To do that, ISOs must consider the cascades topology (considering that upstream water discharges can be used to generate energy in downstream hydros), network transmission constraints (considering that the energy produced by generators must flow through the network to supply demand), thermoelectric generation capacity and fuel costs, and the inflow uncertainty\footnote{{Many other uncertainties can be considered, but for the sake of simplicity and didactic purposes we concentrate on the main format with which this problem is generally studied.}}}. In this setting, a dispatch policy or rule must be estimated and then evaluated through an out-of-sample implementation scheme. In this section, we show both the estimation and evaluation procedures of LDRs applied to LHDP, which is a relevant, representative, and didactic case of a MSLP.

\subsection{General formulation for the MSLP and linear decision rules}

The general form of the MSLP we will approach in this paper follows the one used in \cite{Bodur_2018} and can be written as follows:
    \begin{align}
        &\text{min} \; \mathbb{E}_{\omega}\Bigg[\sum_{t\in \mathcal{T}} \Big(c'_{t}(\xi_{t}(\omega))x_{t}(\xi_{[t]}(\omega)) + q'_{t}(\xi_{t}(\omega))y_{t}(\xi_{[t]}(\omega))\Big)\Bigg] \label{General_ObjFun} \\
        &\text{s.t.:}  \;\big[x_{t}(\xi_{[t]}(\omega)),y_{t}(\xi_{[t]}(\omega))\big] \in   \mathcal{F}_{t}\big(x_{t-1}(\xi_{[t-1]}(\omega)),\xi_{t}(\omega)\big), \quad \forall{t,\omega}. \label{General_Const} 
    \end{align}

In this general model, {$c'_t$ and $q'_t$ stand for the transpose of vectors $c_t$ and $q_t$}, $t$ represents the time steps or decision stages assumed to belong to $\mathcal{T}$ and $\omega$ represents a scenario (or point) of the sample space $\Omega$. The uncertainty is accounted for by random vectors $\xi_{t}:\Omega\longrightarrow \Xi_t$ for all $t\in\mathcal{T}$, whose history until $t$ is represented by vector $\xi_{[t]}=[\xi_{1},...,\xi_{t}]$. {We assume that the support set of $\xi_t$, namely, $\Xi_t$, is a subset} of $\mathbb{R}^{d_{\xi_t}}$ and $\Xi_{[t]}=\bigtimes_{k=1}^t\Xi_k$. In this context, $\mathbb{E}_{\omega}$ means the mathematical expectation, where $\omega$ is used to clearly show where the expectation is acting. Then, for instance, the expected value of $\xi_t$ is defined as $\mu_t = \mathbb{E}_{\omega}[\xi_t(\omega)] = \int_{\omega\in\Omega}\xi_t d\mathbb{P}$. We consider random variables with finite first and second moments. Regarding the decisions, $x_{t}:\Xi_{[t]}\longrightarrow\mathbb{R}^{d_{x_t}}$ is part of the decision vector that stands for the state of the system at the end of stage $t$, whereas $y_{t}:\Xi_{[t]}\longrightarrow\mathbb{R}^{d_{y_t}}$ represents the current stage vector of control actions. Regarding the problem structure, $\mathcal{F}_{t}$ represents a polyhedral set parameterized on the value of the state vector at the end of $t-1$, or initial state of $t$, and on the history of the uncertainty until $t$. Therefore, expression \eqref{General_Const} means a general affine relation among vectors $x_{t}$, $y_{t}$, $x_{t-1}$, and $\xi_{[t]}$.  

Note that, in the general formulation \eqref{General_ObjFun}--\eqref{General_Const}, both the state and the control vectors of a given stage $t$ are general functions of $\xi_{[t]}(\omega)$, the history of the stochastic process until $t$. As per \cite{Bodur_2018}, the LDR approach assumes both the state and the control to belong to a functional space that can be expanded by a basis of nonlinear functions. Alternatively, the two-stage LDR approach assumes that only state variables follow the LDR. In this case, control actions are afterward obtained by a state-target tracking implementation procedure. Finally, it is worth mentioning that although the proposed approach applies to the more general model \eqref{General_ObjFun}--\eqref{General_Const}, in this work, we showcase our contributions based on a particularized version of this problem. We use the LHDP as an application to derive experiments and test the hypothesis that AdaLASSO-regularized LDRs provide improved out-of-sample results in comparison to their non-regularized version.


\subsection{Linear decision rule {applied to the LHDP}}
\label{sec:LHDP}
In this work, we focus on the two-stage LDR approach \cite{Bodur_2018} to generate a hydrothermal scheduling policy considering uncertainty in water inflow. Although many variants of this problem can be used, we keep the model simple (following previously reported works) to study the main properties of an applied LDR in a didactic yet meaningful problem. Therefore, we use the following operational model {to estimate} the LDR:
\begin{align}
    \text{min} &\; \mathbb{E}_{\omega}\Bigg[\sum_{t\in \mathcal{T}}{\alpha^{-t}} \Big(c_t' g_{t}(\omega) + c_{d}'\delta_{t}(\omega) \Big)\Bigg] \label{Ori_ObjFun} \\
    \text{s.t.:}\; \; & Tg_{t}(\omega) + H\rho(u_{t}(\omega)) + Af_{t}(\omega) + \delta_{t}(\omega)
    = d_{t}, \quad \forall{t,\omega} \label{Ori_EnergyBal} \\
    &v_{t}(\omega) - v_{t-1}(\omega) + M(u_{t}(\omega) + s_{t}(\omega)) =\xi_{t}(\omega), \quad \forall{t,\omega}  \label{Ori_WaterBal} \\
    &v_{1}(\omega) = v_{0}, \quad \forall{\omega}  \label{Ori_Vi} \\
    &v_{T}(\omega) \ge v_{f}, \quad \forall{\omega}  \label{Ori_Vf} \\
    &\{v_{t}(\omega), u_{t}(\omega), s_{t}(\omega)\} \in \mathcal{H}_{t}, \quad \forall{t,\omega}  \label{Hydro_Const} \\
    &\{g_{t}(\omega), f_{t}(\omega)\} \in \mathcal{N}_{t}, \quad \forall{t,\omega}  \label{Operative_Const} \\
    &v_{t}(\omega) - \sum\limits_{k\in\mathcal{K}}  \Psi_{k}({\xi}_{[t-\tau:t]}(\omega))\Theta_{t,k} = 0, \quad \forall{t,\omega}. \label{Ori_AffineRule}
\end{align}

First, it is important to mention that throughout this work, decision variables are all on the left-hand side of constraints as well as constants are all on the right-hand side. Thus, in model \eqref{Ori_ObjFun}--\eqref{Ori_AffineRule}, the objective function, \eqref{Ori_ObjFun}, represents the expected value of the net present value at a discount factor $\alpha$ of operating costs within the set of periods $\mathcal{T}$. For the sake of simplicity, the operating cost is represented as a linear function of the thermoelectric generation vector, $g_{t}(\omega)$. Therefore, $c_t$ represents the vector of unitary costs of each thermoelectric unit. Additionally, the objective function is also composed of a second term representing the cost of the energy not supplied, $\delta_{t}(\omega)$, where a unitary cost for the deficit of supply, $c_{d}$, is assumed. 

The nodal energy balance is presented in \eqref{Ori_EnergyBal}. In this expression, total supply, composed of thermoelectric generation, hydro generation, transmission flows, and fictitious generation representing deficit, respectively, should be equal to demand in every node. Thus, $T$ and $H$ allocate thermoelectric and hydro generation to respective buses (nodes of the network) and the incidence matrix $A$, which defines the network topology, associates transmission-line flows to the corresponding nodal equations. Finally, $\rho(\cdot)$ defines the vector of hydro production functions, transforming the amount of water discharged into hydros' turbines into actual power generation. Expression \eqref{Ori_WaterBal} represents the water balance equation (state transition function) for each hydroelectric reservoir, scenario $\omega$, and stage $t$. In this expression, $v_{t}(\omega)$ and $v_{t-1}(\omega)$ represent the vector of water stored at the end of period $t$ and $t-1$ for scenario $\omega$, respectively; $\xi_{t}(\omega)$ represents the vector of incoming water (inflows) in each reservoir for scenario $\omega$; $u_{t}(\omega)$ and $s_{t}(\omega)$ represents the water discharged into the turbines that effectively produce energy and water spillage, respectively, for the period $t$ and scenario $\omega$. Finally, ${M}$ is a matrix that defines rivers' topology in the water balance, accounting for upstream and downstream water release in each reservoir balance equation. The initial condition of the stored water in each reservoir is stated in \eqref{Ori_Vi}. In order to avoid the \emph{end of horizon effect}, Equation \eqref{Ori_Vf} is used to maintain the final reservoir levels greater than or equal to a user-defined limit. Constraint \eqref{Hydro_Const} allows defining operative limits for hydropower plants, e.g., minimum and maximum reservoir storage capacity, minimum water discharge, etc. Similarly, Expression \eqref{Operative_Const} represents the set of technical operational constraints for the generation or transmission assets, e.g., considering minimum and maximum generation, transmission line limits and the Second Kirchhoff Law \cite{Brigatto2017} for each period.

Finally, expression \eqref{Ori_AffineRule} defines a general LDR rule driving the state of the system ($v_t(\omega)$) at each stage and scenario, where $\{\Psi_{k}\}_{k\in\mathcal{K}}$ represents a basis of nonlinear functions. This basis is applied to the inflow vector, ${\xi}_{[t-\tau:t]}(\omega)$, which stacks on ${\xi}_{t}(\omega)$ the vectors of $\tau$-previous observed reservoirs' inflows for each $\omega$. The vector of coefficients $\Theta_{t,k}$ is one of the decisions of our model, which is co-optimized with the rest of the dispatch decisions. It accounts for a limited set of lags, $\mathcal{L}=\{0,...,\tau\}$, and elements, $k \in \mathcal{K}$, in the functional basis. In this work, we use a set of polynomial functions (considering the intercept, i.e., $0 \in \mathcal{K}$) for each reservoir. For the sake of simplicity, we disregard cross-terms effects. Notwithstanding, to consider the influence of other reservoirs $h'\in\mathcal{H}\backslash\{h\}$ in the LDR of reservoir $h$, we also consider the sum of their inflows ($\xi^{(-h)}_{t-l+1}:=\sum_{h'\in\mathcal{H}\backslash\{h\}} \xi_{t-l+1,h'}$) in each LDR expression. It is important to highlight that the consideration of cross-terms{, such as $\xi_{t-l+1,h} \cdot \xi_{t-l'+1,h'}$ for different $l$, $l'$, $h$ and $h'$,} would still be valid under the proposed LDR scheme. Note that although nonlinear, cross-terms accounting for products between the inflow data of different lags and reservoirs can be exogenously pre-calculated and transformed into new variables, e.g., $\xi_{t}^{(l,l',h,h')} = \xi_{t-l+1,h} \cdot \xi_{t-l'+1,h'}$, to be considered in the LDR as customary in linear regression. 

Therefore, {model \eqref{Ori_ObjFun}--\eqref{Ori_AffineRule} defines the optimal }LDR for the hydro reservoir $h$, period $t$, and scenario $\omega$ under the following functional form: 
\begin{equation} \label{Mod_AffineRule}
\begin{split}
    v_{t,h}(\omega) = &\sum_{k\in\mathcal{K}}\sum_{l\in\mathcal{L}_k} \theta_{t,h,k,l}^{(h)} \cdot \left(\frac{\xi_{t-l+1,h}(\omega)-\mu_{t-l+1,h}}{\sigma_{t-l+1,h}}\right)^{k} \\
    + &\sum_{k\in\mathcal{K}_0}\sum_{l\in\mathcal{L}} \theta_{t,h,k,l}^{(-h)} \cdot 
    \left(\frac{\xi^{(-h)}_{t-l+1}(\omega)-\mu^{(-h)}_{t-l+1}}{\sigma^{(-h)}_{t-l+1}}\right)^{k}
    \forall{t,h,\omega}.
\end{split}
\end{equation}
Where, $\mu_{t-l+1,h}$ and $\sigma_{t-l+1,h}$ are the expected value and standard deviation of the inflow of hydro reservoir $h$ at stage $t-l+1$, respectively, $\mu^{(-h)}_{t-l+1}$ and $\sigma^{(-h)}_{t-l+1}$ represent the same statistics for $\xi^{(-h)}_{t-l+1}$, and {$\theta_{t,h,k,l}^{(h)}$ are the components of the coefficient vector $\Theta_{t,k}$}. By considering the intercept, $k=0$, seasonal levels are automatically addressed.  Thus, we can use the deseasonalized series to determine the policy based on standardized shocks above or below average. Additionally, we avoid degenerated models (with multiple coefficients playing the role of intercept) by disregarding multiple lags when $k=0$ in the first term, i.e., by making $\mathcal{L}_k:=\mathcal{L} \; \forall k>0$ and $\mathcal{L}_0:=\{0\}$, and disregarding $k=0$ in the second term, i.e., by making $\mathcal{K}_0:=\mathcal{K}\backslash \{0\}$.

\subsection{{Linear decision rule estimation}}
\label{sec:LDRestimation}

In practice, according to \cite{Bodur_2018}, the so-called SAA version of problem \eqref{Ori_ObjFun}--\eqref{Ori_AffineRule} is used to estimate the vector of coefficients $\hat{\Theta}$ for the LDR. To do that, we randomly generate (based on a Monte Carlo procedure) a finite sample space $\Omega_N$ from $\Omega$ with $N$ scenarios, each of which associated with a probability equal to $1/N$.\footnote{{Notwithstanding, scenario reduction techniques, clusterization methods, and variants of Monte Carlo approaches can also be used to generate the aforementioned set of scenarios used to estimate LDRs. Depending on the method, one just needs to consider an individual probability different than $1/N$ for each scenario. Henceforth, however, for the sake of simplicity, we will consider equally distributed scenarios as per the SAA approach.}} We denote as $\{\xi_{t,s}\}_{s\in\Omega_N, t\in\mathcal{T}}$ the associated inflow values for each sampled scenario. Thus, we can estimate the LDR coefficients by solving the following two-stage linear optimization problem:

\begin{align}
    \text{min} &\frac{1}{N}\sum_{s \in \Omega_N}\sum_{t\in \mathcal{T}}{\alpha^{-t}} \Big(c_t' g_{t,s} + c_{d}'\delta_{t,s} \Big) \label{SAAOri_ObjFun} \\
    \text{s.t.:}\; \; & Tg_{t,s} + H\rho(u_{t,s}) + Af_{t,s} + \delta_{t,s}
    = d_{t}, \quad \forall{t,s} \label{SAAOri_EnergyBal} \\
    &v_{t,s} - v_{t-1,s} + M(u_{t,s} + s_{t,s}) = \xi_{t,s}, \quad \forall{t,s}  \label{SAAOri_WaterBal} \\
    &v_{1,s} = v_{0}, \quad \forall{s}  \label{SAAOri_Vi} \\
    &v_{T,s} \ge v_{f}, \quad \forall{s}  \label{SAAOri_Vf} \\
    &\{v_{t,s}, u_{t,s}, s_{t,s}\} \in \mathcal{H}_{t}, \quad \forall{t,s}  \label{SAAHydro_Const} \\
    &\{g_{t,s}, f_{t,s}\} \in \mathcal{N}_{t}, \quad \forall{t,s}  \label{SAAOperative_Const} \\
    &v_{t,s} - \sum\limits_{k\in\mathcal{K}}  \Psi_{k}({\xi}_{[t-\tau:t],s})\Theta_{t,k} = 0, \quad \forall{t,s}. \label{SAAOri_AffineRule}
\end{align}

It is relevant to note that the number of coefficients in the decision rule depends on the number of periods ($|\mathcal{T}|$), water reservoirs ($|\mathcal{H}|$), time-lags of the historical data ($|\mathcal{L}|$), and on the degree of the polynomial basis ($|\mathcal{K}|$). Although it provides great flexibility for obtaining reasonable dispatch policies, the large number of coefficients easily approaches the typical number of scenarios that can be handled in practical applications.\footnote{For instance, \cite{Bodur_2018}, and more recently \cite{daryalal2020lagrangian}, apply LDRs with relatively few scenarios in-sample. For instance, in \cite{daryalal2020lagrangian}, \textit{``we approximated this problem with just 25 scenarios due to computational limitations."} and in \cite{Bodur_2018}, authors use only \textit{``250 scenarios to construct an SAA"}.} For instance, if we consider 24 months (2 years horizon on monthly basis), five reservoirs, 12 lags, and six degrees in the polynomial basis, the total number of coefficients for estimating the LDR is approximately equal to 8.6 thousand. This number (approximately) doubles if we account for the second term of expression \eqref{Mod_AffineRule} (in this case, approximately equals 16.8 thousand) and can get much larger if we consider cross-terms. \emph{As a result, in practical applications where a SAA of the model is used with, say, hundreds to a few thousands of sampled scenarios, it is very likely spurious coefficients will be found due to overfit.} The consequence, in general, is a poor performance when implementing the policy in out-of-sample scenarios as largely reported in regression and machine learning literature.

\subsection{Linear decision rule implementation and out-of-sample evaluation}
\label{sec:outofsample}

The implementation of a two-stage LDR can follow a \textit{state-target tracking} (STT) procedure described in \cite{Bodur_2018}. In this procedure, the operation of the system is implemented by a STT optimization problem, which aims to follow the targeted storage defined by the estimated LDR as much as possible for a given initial state and observed data. The estimated LDR is represented by the optimal coefficient vector $\hat{\Theta}$, which can be obtained as the solution of problem \eqref{SAAOri_ObjFun}--\eqref{SAAOri_AffineRule}. Thus, given the estimated LDR coefficient vector $\hat{\Theta}$, the vector of initial storage $v_{t-1,s}$, and the observed scenario of inflows $\xi_{[t-\tau:t],s}$, the operation of the system for stage $t$ can be defined by the solution of the following deterministic and single-period STT problem:
\begin{align}
    \text{min} & \;\; c' g_{t,s} + c_{d}'\delta_{t,s} + \gamma'(e^{+}_{t,s} + e^{-}_{t,s}) \label{STT_Ori_ObjFun} \\
    \text{s.t.:}\; \; & Tg_{t,s} + H\rho(u_{t,s}) + Af_{t,s} + \delta_{t,s}
    = d_{t} \label{STTOri_EnergyBal} \\
    &v_{t,s} + M(u_{t,s} + s_{t,s}) = v_{t-1,s} + \xi_{t,s}  \label{STT_Ori_WaterBal} \\
    &\{v_{t,s}, u_{t,s}, s_{t,s}\} \in \mathcal{H}_{t} \label{STT_Hydro_Const} \\
    &\{g_{t,s}, f_{t,s}\} \in \mathcal{N}_{t}  \label{Imp_Operative_Const} \\
    &e^{+}_{t,s}, \; e^{-}_{t,s} \ge 0 \\
    &v_{t,s} + e^{+}_{t,s} - e^{-}_{t,s} = \sum\limits_{k\in\mathcal{K}}  \Psi_{k}({\xi}_{[t-\tau:t,s]})\hat{\Theta}_{t,k}. \label{Imp_Ori_AffineRule}
\end{align}
Problem \eqref{STT_Ori_ObjFun}--\eqref{Imp_Ori_AffineRule} delivers as output optimal implemented operational decisions $(g_{t,s}^*,u_{t,s}^*,f_{t,s}^*,\delta_{t,s}^*)$ and the new state of the system $v_{t,s}^*$. It is worth highlighting that $v_{t,s}^*$ can differ from the targeted state defined by the LDR in the right-hand-side of expression \eqref{Imp_Ori_AffineRule}. However, aiming to enforce the LDR as much as possible in the implementation step, a high penalty vector $\gamma$ is considered in the objective function to penalize the absolute value of the state-target deviation errors $e^{+}_{t,s} + e^{-}_{t,s}$. Low penalties for spillages can also be considered for practical purposes albeit omitted in this model. It is relevant to emphasize at this point that the above procedure represents the actual two-stage decision rule. In other words, it provides the system operator with a process that, given $\hat{\Theta}$, the initial state, and the observed data until stage $t$, recommends what to do in this stage in a nonanticipative fashion.

Thus, based on the STT problem \eqref{STT_Ori_ObjFun}--\eqref{Imp_Ori_AffineRule}, we can evaluate out of sample the performances of any LDR for a larger and new set of randomly generated scenarios $\Omega_M$ {(generated from the same underlying process)}, where $M >> N$. To do that, we initialize $v_{0,s}^*$ with a known initial storage vector $v_0$ and successively solve the STT problem for all $t\in\mathcal{T}$ and $s\in\Omega_M$. In each step, i.e., for each $t,s$, we consider as input $(v_{t-1,s}, \xi_{[t-\tau]:t],s}, \hat{\Theta})$ and as output $(g_{t,s}^*,u_{t,s}^*,f_{t,s}^*,\delta_{t,s}^*, v_{t,s}^*)$. Note that the storage output of a given period is used as input for the next period. At the end of this chronological simulation process, we can evaluate the out-of-sample operational cost due to the implementation of a LDR defined by $\hat{\Theta}$ as follows:
\begin{align}
    z^*_M(\hat{\Theta}) = \frac{1}{M}\sum_{s \in \Omega_M}\sum_{t\in \mathcal{T}}{\alpha^{-t}} \Big(c' g_{t,s}^* + c_{d}'\delta_{t,s}^* \Big). \label{CustoOut}
\end{align}
%
Because {the STT is always limited for any value of $\hat{\Theta}$, due to the uniform law of large numbers, under mild conditions on the distribution of $\xi$, we have that $z^*_M(\hat{\Theta}) \longrightarrow \mathbb{E}_{\omega}\Big[\sum_{t\in \mathcal{T}}{\alpha^{-t}} (c' g_{t}^*(\omega) + c_{d}'\delta_{t}^*(\omega))\Big]$ when $M \rightarrow \infty$, in which $g_{t}^*(\omega)$, $\delta_{t}^*(\omega)$, and $v_{t}^*(\omega)$ are the optimal solutions of the STT problem for $\xi_t(\omega)$ (\cite{saa99}).}

\section{The Proposed Regularized Linear Decision Rule via AdaLASSO} \label{Regularized_AffinePol}

The AdaLASSO is a regularization operator that has the oracle property \cite{zou2006adaptive}. Roughly, it means that such a regularization scheme is capable of selecting only a few coefficients truly relevant for explaining the data. This concept is proposed for the first time in the present work to improve the quality of LDR applied to solve MSLP problem. The proposed regularized LDR estimation model considers an additional term in the objective function penalizing the scaled $\ell_1$-norm of the coefficient vector. Then, the new objective function accounts for the trade-off between the in-sample cost minimization and the number of parameters (parsimony) used in the LDR to minimize it. The goal of the regularization is to obtain a better generalization of the LDR for unseen scenarios, thereby reducing the out-of-sample cost.

If $\Theta$ stacks all the coefficients $\theta_{t,h,k,l,r=1}=\theta_{t,h,k,l}^{(h)}$ and $\theta_{t,h,k,l,r=2}=\;\theta_{t,h,k,l}^{(-h)}$ for all $t\in\mathcal{T}$, $h\in\mathcal{H}$, $k\in\mathcal{K}$, and $l\in \mathcal{L}$, the AdaLASSO regularization function $F_{\lambda,\Theta^{(0)}}(\Theta)$ can be defined as follows:
\begin{equation} \label{PenaltyFunction}
    F_{\lambda,\Theta^{(0)}}(\Theta) = \lambda\sum_{t\in\mathcal{T}}\sum_{h\in\mathcal{H}}\sum_{k\in\mathcal{K}_0}\sum_{l\in\mathcal{L}}\sum_{r=1}^2 \frac{|\theta_{t,h,k,l,r}|}{|\theta_{t,h,k,l,r}^{(0)}|}. 
\end{equation}
Function \eqref{PenaltyFunction} receives a vector of parameters $\Theta$ and returns a scaled sum of the absolute values of its components, disregarding the intercept coefficient ($k=0$). It has two parameters, namely: 1) $\lambda$, which is a scalar reflecting the overall penalization level to all coefficients, and 2) $\Theta^{(0)}$, which { is the vector that stacks the non-regularized LDR coefficients $\theta_{t,h,k,l,r}^{(0)}$ estimated with the original SAA problem \eqref{SAAOri_ObjFun}--\eqref{SAAOri_AffineRule}, where zero elements are replaced with 1}. 
It is important to notice that from one period to another, the value of $\Theta^{(0)}$ should not dramatically change. Therefore, in practical applications, one might not need to solve the two problems all the time. {Additionally, it is important to highlight that other interesting regularization metrics can also be used as alternatives to the AdaLASSO in cases the extra computational burden to estimate $\Theta^{(0)}$ becomes an issue.}\footnote{{Note that many other regularization penalties could be used, each of which with their pros and cons as widely reported in the related literature. For instance, a combination of $\ell_1$ and $\ell_2$ norms (Elastic Net) is a salient example that could be used as an alternative. Notwithstanding, it is beyond the scope of this paper to analyze all of them. Instead, we focus on the more general weighted $\ell_1$-norm form of the AdaLASSO, which is widely used and recognized as a relevant tool for variable selection \cite{medeiros2020} and to improve the model performance in out-of-sample tests.}} Nevertheless, as we are in a planning context, the computational time to solve both models might not constitute an issue. 

Thus, the (AdaLASSO) regularized LDR for the MSLP can be estimated by the solution of the following two-stage linear programming model: 
\begin{align}
    \text{min} & \;\displaystyle\frac{1}{N}\sum\limits_{s\in\Omega_N}\sum\limits_{t\in\mathcal{T}}{\alpha^{-t}}\Big( c' g_{t,s} + c'_{d}\delta_{t,s}\Big) + \lambda\sum_{t\in\mathcal{T}}\sum_{h\in\mathcal{H}}\sum_{k\in\mathcal{K}_0}\sum_{l\in\mathcal{L}}\sum_{r=1}^2 \frac{\phi_{t,h,k,l,r}}{|\theta_{t,h,k,l,r}^{(0)}|} \label{Mod_ObjFun} \\
    &\text{s.t.: Constraints \eqref{SAAOri_EnergyBal} to \eqref{SAAOri_AffineRule}} \\
    &\phi_{t,h,k,l,r} - \theta_{t,h,k,l,r} \geq 0, \qquad \forall{t,h,k,l,r} \label{Mod_ThetaPos} \\
    &\phi_{t,h,k,l,r} + \theta_{t,h,k,l,r} \geq 0, \qquad \forall{t,h,k,l,r} \label{Mod_ThetaNeg} 
\end{align}
In \eqref{Mod_ObjFun}--\eqref{Mod_ThetaNeg}, the regularization term defined in \eqref{PenaltyFunction} is added to the original objective function, \eqref{SAAOri_ObjFun}. Additionally, the nonlinear absolute-value terms, $|\theta_{t,h,k,l,r}|$, in \eqref{PenaltyFunction} are replaced with auxiliary variables, $\phi_{t,h,k,l,r}$, and epigraph constraints, \eqref{Mod_ThetaPos}--\eqref{Mod_ThetaNeg}, following the standard linearization procedure in linear programming. {Note that the output of \eqref{Mod_ObjFun}--\eqref{Mod_ThetaNeg} is a vector of coefficients $\hat{\Theta}_\lambda$ representing a regularized LDR, which can be evaluated out of sample through \eqref{CustoOut} for different values of $\lambda$. In this context, a grid search based on the out-of-sample cost metric can be used to obtain the best regularized LDR within all policies induced by the inspected values of $\lambda$.} Finally, it is important to highlight that although we have used the LHDP problem to derive our ideas, the results of this section are very general and can be easily applied to other MSLP instances particularized from \eqref{General_ObjFun}--\eqref{General_Const}.

\section{Computational Experiments} \label{Quantitative_Analysis}

For analyzing the proposed method, two different case studies are performed. The goals of the first case study are 1) to study the quality of the regularized LDR in out-of-sample tests, and 2) to study the properties of regularized coefficients. To address the former, we benchmark the regularized LDR with the classical non-regularized LDR \cite{Bodur_2018} and evaluate the out-of-sample costs. Finally, the latter is addressed based on comparisons between the coefficients obtained with the non-regularized and regularized LDRs. The second case study shows statistical analysis of the out-of-sample results of the regularized LDR based on a more complex cascade topology for a stochastic process with nonlinear temporal dependencies \cite{GAS}. It focuses on the effects of the regularization process over the total cost distribution, for which we show further benefits obtained from the proposed regularization beyond cost reduction. 

All computational experiments were conducted using the Julia programming language (version v1.2) and solved with Gurobi v0.9.13, on an Intel Core i7-10510U processor at 1.80 GHz with 8 GB of RAM. 

\subsection{Case Study 1 -- a simple test system} \label{Benchmark_Comparison}
To evaluate the proposed regularized policy based on LDR and its performance in reducing the overfitting issue, we consider a comparison between the proposed method and the standard non-regularized LDR approach. For the sake of simplicity and comparison purposes, in this case study, we use a stagewise-independent stochastic process with periodic mean and variance based on historical data to generate inflow scenarios. Thus, in this particular case, the standard SDDP approach provides the optimal policy, thereby serving as a reference for the proposed LDR approach. To this end, we implement the same operational model used to define each stage of model \eqref{Ori_ObjFun}--\eqref{Ori_AffineRule} in \cite{dowson_sddp.jl} and used the same in-sample scenarios.

In this first simulation, a simple system configuration {, with only one reservoir (one state variable),} is chosen to ensure we can explore and study operational results obtained with the regularized LDR. Table \ref{tab:case1_system_config} presents the main characteristics of the simulation.

\begin{table}[H]
    \centering
    \caption{Case Study 1 - System Configuration}
    \label{tab:case1_system_config}
    \begin{tabular}{c|c}
         \hline
         {Number of Stages}        & 36    \\
         \hline
         {Number of Hydros}        & 1     \\
         \hline
         {Number of Thermals}      & 6     \\
         \hline
         {In-sample Scenarios}     & 100   \\
         \hline
         {Out-of-sample Scenarios} & 1000   \\
         \hline
         {Discount Rate} & 0.5\%          \\
         \hline
    \end{tabular}
\end{table}

The main characteristics of the hydropower plant considered in this case study are presented in Table \ref{tab:case1_hydro_param}, while thermal projects are presented in Table \ref{tab:case1_thermal_param}.

\begin{table} [H]
    \centering
    \caption{Case Study 1 - Hydro Power Plants Characteristics}
    \label{tab:case1_hydro_param}

     \begin{tabular}{c|c}
         \hline
         {Max Reservoir Storage (units)}  & 894    \\
         \hline
         {Min Reservoir Storage (units)}  & 356     \\
         \hline
         {Max Turbine Capacity (units)}  & 700     \\
         \hline
         {Efficiency (MW/unit)}          & 0.414   \\
         \hline
    \end{tabular}
\end{table}

\begin{table} [H]
    \centering
    \caption{Case Study 1 - Thermal Power Plants Characteristics}
    \label{tab:case1_thermal_param}
    \begin{tabular}{c|c|c}
         \hline
         {Thermal} & {Capacity (MW)} & {Variable Cost (\$/MWh)} \\
         \hline
         {Thermal 1}  & 100 & 500   \\
         {Thermal 2}  & 100 & 113   \\
         {Thermal 3}  & 100 & 153    \\
         {Thermal 4}  & 50  & 116      \\
         {Thermal 5}  & 50  & 58    \\
         {Thermal 6}  & 50  & 86   \\
         \hline
    \end{tabular}
\end{table}

In this case, where only one hydro unit is considered, the second term in \eqref{Mod_AffineRule} is not applicable. Thus, the parameters of the LDR will be selected considering the maximum {order in the polynomial basis equal to six (i.e., $\max\{\mathcal{K}\}=6$) and the maximum value of lag is 11 (i.e., $\max\{\mathcal{L}\}=11$)}.

\subsubsection{Policy quality}

To measure the quality of the policies generated with the regularized LDR, we use the optimization model proposed in \eqref{Mod_ObjFun}--\eqref{Mod_ThetaNeg} for different values of $\lambda$. The cost of each policy is evaluated considering out-of-sample scenarios according to the STT procedure described in Section \ref{sec:outofsample}. The best regularized LDR is selected according to the minimum out-of-sample cost metric \eqref{CustoOut} {by varying the value of $\lambda$ within a user-defined grid}\footnote{If more than one value of $\lambda$ in the grid exhibit the same out-of-sample cost value, we arbitrarily select the first one, with the lowest regularization penalty.}. Table \ref{tab:ComputationalTime_Case1A} shows the computational times required in the estimation and out-of-sample evaluation processes as well as the respective expected costs.

\begin{table}[H]
    \centering
    \caption{Computational Time - Case Study 1}
    \label{tab:ComputationalTime_Case1A}
    \begin{tabular}{c|c|c|c|c}
        \hline
         & \multicolumn{2}{c|}{Comp. Time} & \multicolumn{2}{c}{Exp. Total Cost}\\
         & Estimation & Evaluation & Estimation & Evaluation\\
& (s) & (s) & (M\$) & (M\$)\\
        \hline
        {$\lambda = 0$} & 9.73 & 5.75 & 1.56 & 2.00 \\
        \hline
        {$\lambda = 10^3$} & 4.60 & 8.02  & 1.59 & 1.63 \\
        \hline
    \end{tabular}
\end{table}

\begin{figure}[H]
    \centering
    \includegraphics[scale=.5]{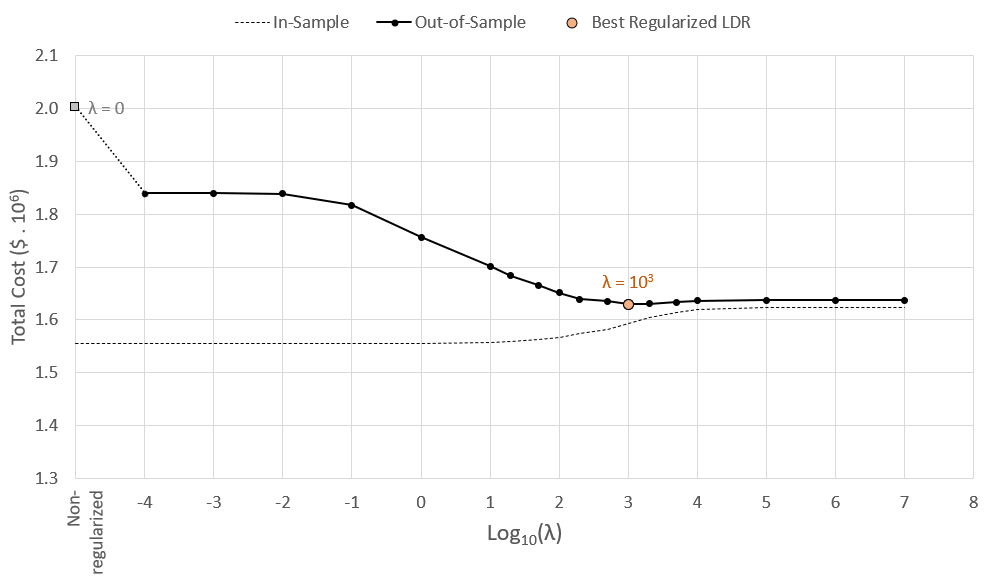}
    \vspace{-0.5cm}
    \caption{Total in-sample and out-of-sample system cost per $\lambda$}
    \label{fig:Case1_TotalCost_Average}
\end{figure}

Figure \ref{fig:Case1_TotalCost_Average} shows the in-sample (estimation) and out-of-sample (evaluation) system operating cost for different values of $\lambda$. 
Noticeably, the in-sample operational cost does not significantly vary with values of $\lambda$ {as the out-of-sample does}. Additionally, the AdaLASSO penalization, albeit shrinking parameters, do not significantly affect the {in-sample cost} until a determined threshold. In turn, as the value of $\lambda$ increases, the out-of-sample operating cost significantly decreases{, revealing the existence of many regularized policies that can significantly outperform the traditional non-regularized LDR.} These facts show us evidence of overfitting and the potential threat of non-regularized LDR in finite samples. Furthermore, it also shows that, in finite samples, better policies, with improved out-of-samples results, may be found based on the proposed regularization scheme.  

The increase of the regularization penalty shrinks a great number of parameters associated with \emph{explanatory variables} (using the terminology of regression analysis in statistics) that are not indeed relevant to generalize the LDR model to other samples. For the defined system configuration (classified loosely restricted case -- we will further show a more restricted case) and stage-wise independent scenarios, the lower out-of-sample cost is found at $\lambda = 10^3$. However, the out-of-sample results for $\lambda > 10^3$ do not vary significantly until the maximum $\lambda$ (which would mean that the LDR considers only the intercepts). In this case, the proposed regularization scheme provides an $18.5\%$ gain, in terms of out-of-sample cost reduction, compared to the policy obtained with the traditional non-regularized LDR approach. This gain grows up to 25\% for the 95-percentile of the cost distribution.

Regarding the benchmark with the SDDP reference, we use the same in-sample scenarios used to train the LDRs as the set of backward scenarios in the SDDP. Forward steps are calculated using random samples also from the set of in-sample scenarios. After the SDDP algorithm converged (within a computational time of 99.9 seconds), we use the same 1000 out-of-sample scenarios used to evaluate the LDR approaches and obtained a reference cost of 1.57 M\$. While the in-sample cost obtained with the non-regularized LDR falls within the in-sample GAP of the SDDP, the out-of-sample cost {of the non-regularized LDR is 27.39 \% higher than the out-of-sample of the SDDP}. On the other hand, the best regularized LDR exhibits an improved out-of-sample performance, { only $3.8\%$} higher than the SDDP reference {(which is 7.16 times smaller than the gap obtained with the non-regularized LDR)}.



\subsubsection{Operational results}

\begin{figure}[h!]
    \centering
    \includegraphics[scale=.5]{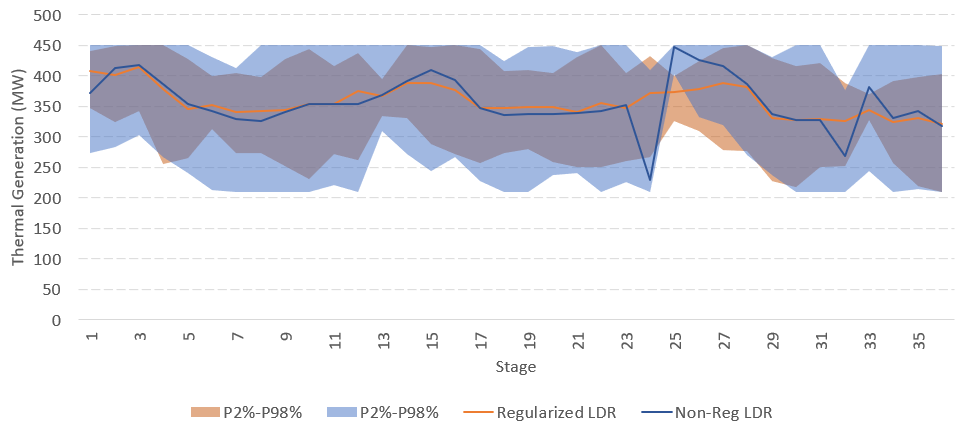}
    \caption{Thermal power generation in avgMW}
    \label{fig:Case1_GerTer}
\end{figure}

\begin{figure}[h!]
    \centering
    \includegraphics[scale=.5]{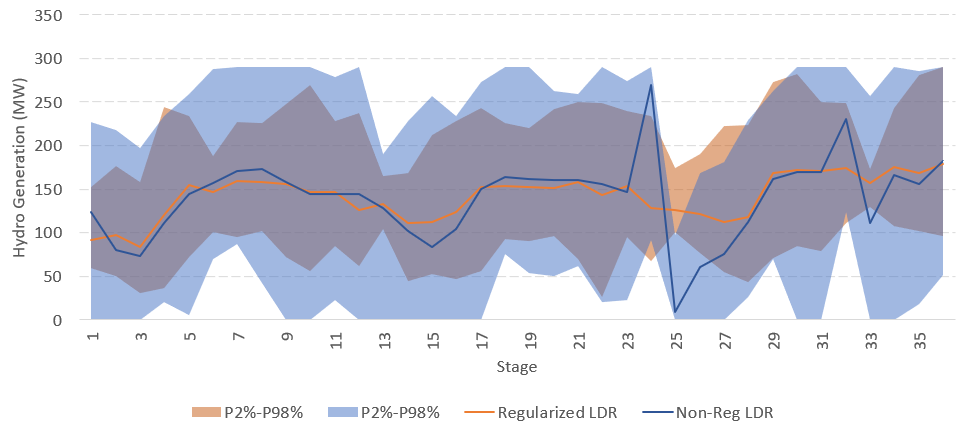}
    \caption{Hydro power generation in avgMW}
    \label{fig:Case1_GerHid}
\end{figure}

\begin{figure}[h!]
    \centering
    \includegraphics[scale=.5]{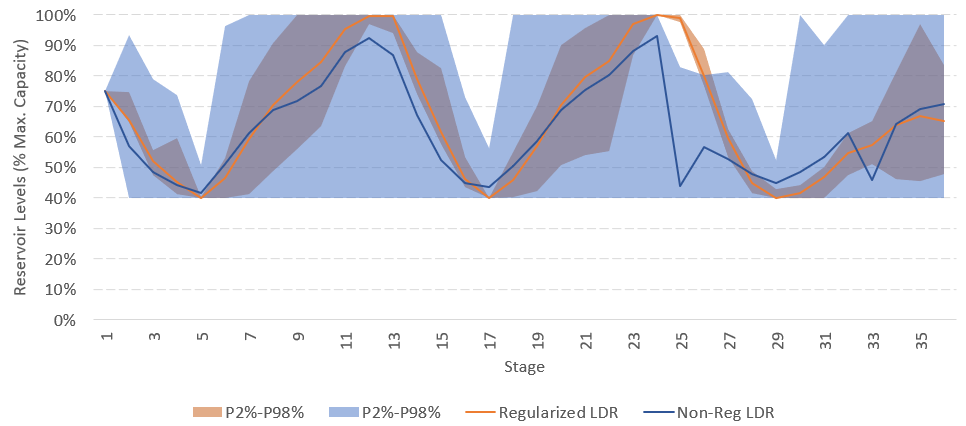}
    \caption{Hydropower plant initial level of the reservoirs at each stage in percentage of the maximum}
    \label{fig:Case1_Volume}
\end{figure}

In this subsection, we compare operative results simulated with the non-regularized (lines and shadows in blue) and regularized (lines and shadows in orange) LDR. Hereinafter, whenever we use the term regularized LDR, we mean the best regularized LDR, i.e., with the best value of $\lambda$ previously selected (in this case, $\lambda=10^3$). Figures \ref{fig:Case1_GerTer}, \ref{fig:Case1_GerHid}, and \ref{fig:Case1_Volume} show the total thermal generation, hydro generation, and hydro reservoir storage levels, respectively, for the out-of-sample scenarios. Lines represent expected values, whereas shadows account for the 2\% to 98\% quantile range. 

Results obtained from both policies exhibit similar average levels, trends, and seasonalities. Nevertheless, moderate deviations in the average amplitude of the seasonal pattern reveal differences in the overall storage strategy, whereas prominent spikes in some stages highlight the instability that non-regularized coefficients may bring under finite sample sizes. Our tests indicate that small regularization penalties can mitigate such cases, suggesting that they are mostly caused by the combination of overfitting and degenerated solutions. In addition to the aforementioned deviations, the most prominent difference between the two policies lies in their variability. Noticeably, the non-regularized policy exhibits a much more uncertain profile (see blue shadows in Figures \ref{fig:Case1_GerTer}, \ref{fig:Case1_GerHid}, and \ref{fig:Case1_Volume}) if compared with the proposed regularized policy. In some stages, where the storage levels approach the maximum and minimum bounds, only the intercept of the LDR is selected, thereby resulting in a deterministic policy in which storage targets are defined regardless of the inflow. On the other hand, the non-regularized dependency on inflows under unseen scenarios can produce expensive thermal generation. For instance, in Figure \ref{fig:Case1_GerTer}, the activation of the most expensive thermal resource induced by the non-regularized LDR imposes a cost many times higher than the cost savings obtained with lower generation scenarios (see the generation costs in Table \ref{tab:case1_thermal_param}).



\subsubsection{Regularized coefficients}

\begin{figure}[H]
    \centering
    \includegraphics[scale=.5]{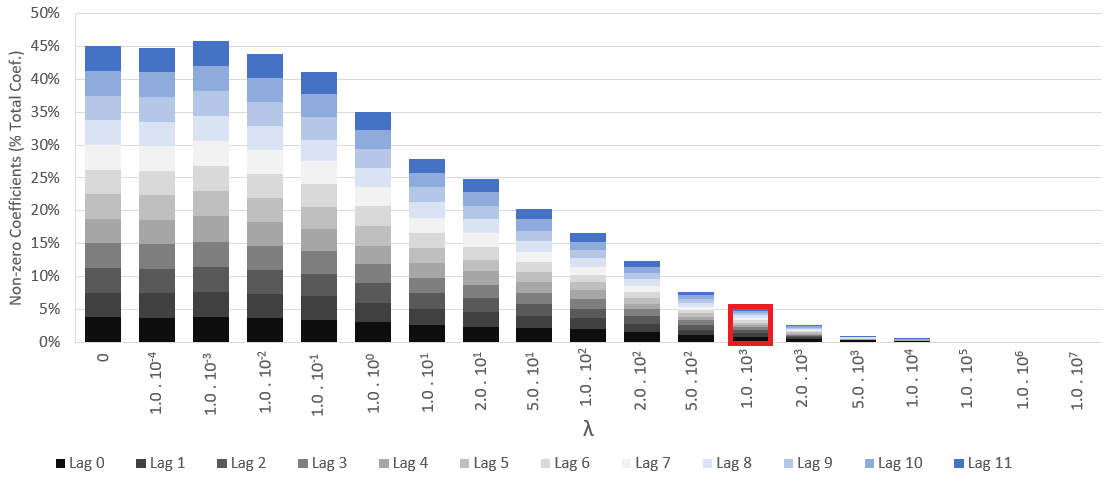}
    \caption{Coefficient selection in percentage of the total available coefficients}
    \label{fig:Case1_CoefSelection}
\end{figure}

Finally, we provide an analysis of the selection and shrinkage of the decision variables (LDR parameter). Figure \ref{fig:Case1_CoefSelection} illustrates the percentage of non-zero coefficients defined for each (policy) value of $\lambda$ per lag of the water inflow. Noticeably, the non-regularized policy ($\lambda = 0$) contains a relatively large number of coefficients different than zero, e.g., $45.0\%$ of the total number of coefficients. Increasing the penalty value, however, the number of non-zero coefficients significantly diminishes, even though the in-sample total cost remains approximately the same (as can be observed in Figure \ref{fig:Case1_TotalCost_Average}). This reduction in the number of LDR coefficients clearly indicates the existence of degenerated solutions -- solutions with the same in-sample operational cost. 

Additionally, the existence of different solutions with the same in-sample cost but with a reduced number of parameters and cheaper out-of-sample costs strongly corroborates the main hypotheses of this work\footnote{As per the introduction section, the main hypothesis is: \emph{non-regularized LDR overfits in-sample data and more parsimonious LDR-based policies estimated with AdaLASSO-regularization can provide relevant cost savings under unseen out-of-sample data.}}. At the best penalty $\lambda=10^3$ (highlighted in red in Figure \ref{fig:Case1_CoefSelection}), only 5.1\% coefficients are different than zero, which represents a relevant reduction in comparison to the 45.0\% obtained with $\lambda = 0$ (i.e., an 88.6\% reduction in the non-zero coefficients).

\begin{figure}[H]
    \centering
    \includegraphics[scale=.5]{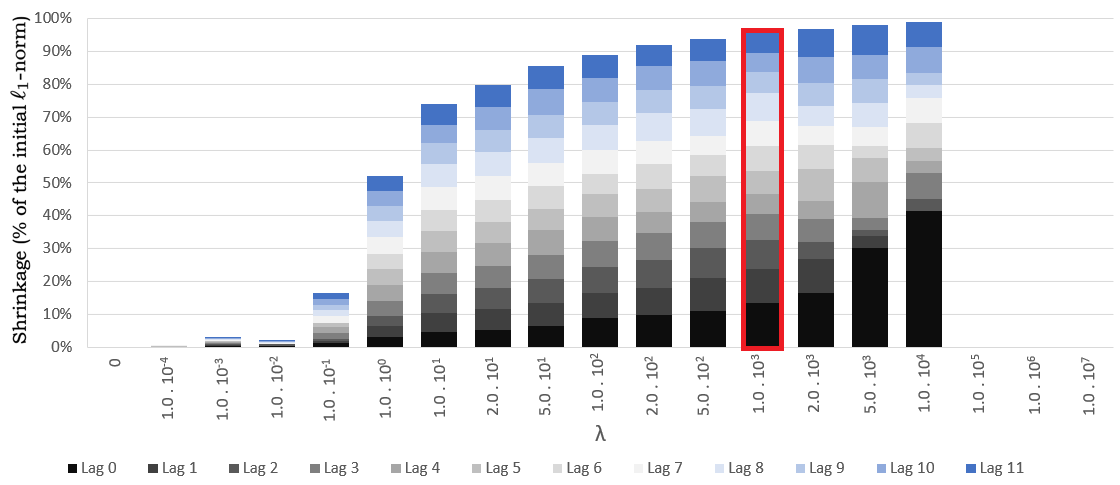}
    \caption{Coefficient $\ell_1$-norm shrinkage in percentage of the $\ell_1$-norm for $\lambda = 0$}
    \label{fig:Case1_CoefShrinkage}
\end{figure}

Figure \ref{fig:Case1_CoefShrinkage} depicts the shrinkage process in the LDR coefficients. It is interesting to see that the $\ell_1$-norm of vector $\Theta^{(\lambda=10^3)}$ is 95.57\% lower than the same norm calculated with $\lambda = 0$.

\subsection{Case Study 2}

For the sake of evaluating the performance of the regularized LDR in a more realistic problem, Case Study 2 is carried out considering the operation of interconnected cascades. The water inflows into the reservoirs are modeled via generalized autoregressive with score model (GAS), a non-Gaussian stochastic process with nonlinear score-driven autoregression dependence. This framework allows the modeling and simulation of inflow time series within a coherent framework that, e.g., recognizes the boundaries and specific family of conditional distributions of the underlying stochastic process. We refer the interested reader to \cite{henrique} for an example of an application in wind power generation in Brazil and to \cite{GAS} for the open-source tool used to generate the scenarios for this case study. The consideration of realistic processes with typical nonlinear dependencies and non-Gaussian distributions constitutes a relevant and timely avenue for approximating the mathematical model and the true underlying problem. It is worth mentioning that the usage of a stochastic process with non-linear (non-convex) time dependencies prevents the application of the traditional SDDP technique. {Interestingly, the LDR considers the stochastic processes driving the uncertainties as a completely exogenous model. In this case, all sorts of spatial-temporal dependencies can be considered through simulated scenarios and used to estimate and evaluate the LDR as per previous sections.}

In this case study, we further analyze the economic benefits of the proposed idea of regularizing LDRs. Therefore, we extend the analysis previously carried out for the expected cost to analyze the cost distribution. Finally, we also analyze specific operational and market performance indices to exemplify other economic benefits of regularizing LDRs.

\subsubsection{System Configuration}

\begin{figure}[h!]
    \centering
    \includegraphics[scale=.6]{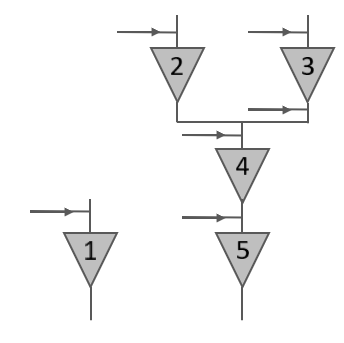}
    \caption{Hydropower plants cascade topology}
    \label{fig:Case2_Hydroplants_Config}
\end{figure}

In this Case Study, the system configuration is composed of five reservoirs {(five state variables)} in two different cascades as shown in Figure \ref{fig:Case2_Hydroplants_Config}. The main hydro power plants' data are depicted in Table \ref{tab:Case2_Hydro_Config}.

\begin{table}[h!] 
    \centering
    \caption{Hydropower plants configuration}
    \begin{tabular}{c|c|c|c|c|c}
    \hline
         {Hydro} & {1} & {2} & {3} & {4} & {5} \\
    \hline
         {Max. Storage (hm3)} & 394 & 319 & 291 & 197 & 166 \cr
         {Min. Storage (hm3)} & 36 & 78 & 17 & 21 & 8 \cr
         {Max. Turb. (m3/s)} & 80 & 103 & 77 & 227 & 277 \cr
         {Prod. Factor (MW/month.m3.s)} & 0.81 & 1.12 & 1.10 & 1.10 & 1.19 \cr
         {Downstream Plant} & None & 4 & 4 & 5 & None \cr
         \hline
    \end{tabular}
    \label{tab:Case2_Hydro_Config}
\end{table}

The demand is set as $1000MW$ in every stage, preventing any extra effect on the results besides the inflow scenarios and decisions on the power system scheduling. The load shedding cost is set to $\$2611/MWh$. Moreover, the system configuration of this example contains six thermal power units whose data are shown in Table \ref{tab:Case2_Therm_Config}.

\begin{table}[h!] 
    \centering
    \caption{Thermal power plants configuration}
    \begin{tabular}{c|c|c}
    \hline
         {Name} & {Max. Gen. Capacity} & {Variable Cost} \\
           & (MW) & (\$/MWh) \\ 
    \hline
         Thermal 1 & 250 & 159    \cr
         Thermal 2 & 50 &  113   \cr
         Thermal 3 & 250 & 153   \cr
         Thermal 4 & 50 &  116  \cr
         Thermal 5 & 50 &  58  \cr
         Thermal 6 & 50 &  86  \cr
         \hline
    \end{tabular}
    \label{tab:Case2_Therm_Config}
\end{table}

Finally, Table \ref{tab:Case2_Table_Config} depicts other relevant parameters used in this case study. 

\begin{table}[h!]
    \centering
    \caption{Case study main configurations.}
    \begin{tabular}{c|c|c|c}
    \hline
    {Stages} & {In-Sample Series} & {Out-of-Sample Series} & {Max. Lag} \\
    \hline
    24 & 100 & {10000} & 12 \\
    \hline
    \end{tabular}
    \label{tab:Case2_Table_Config}
\end{table}

\subsubsection{Regularization Process}

As proposed, the regularized LDR method requires the definition of an optimal penalty for the AdaLASSO, therefore the regularization process was carried out. A similar sequential process of changing the penalty level in the estimation step and evaluating the LDR policy for out-of-sample scenarios was performed. Table \ref{tab:ComputationalTime_Case2} depicts the computational times required in the estimation and out-of-sample evaluation processes for both the non-regularized and best regularized LDRs. 

\begin{table}[H]
    \centering
    \caption{Computational Time - Case Study 2}
    \label{tab:ComputationalTime_Case2}
    \begin{tabular}{c|c|c}
        \hline
         & Estimation Time (s) & Evaluation Time (s)\\
        \hline
        {$\lambda = 0$} & {138.11} & {113.12} \\
        \hline
        {$\lambda = 10^4$} & {29.77} & {79.82} \\
        \hline
    \end{tabular}
\end{table}

\begin{figure}[h]
    \centering
    \includegraphics[scale=.8]{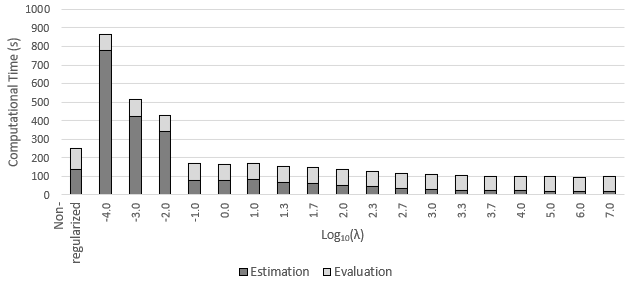}
    \caption{Computational time for each $\lambda$}
    \label{fig:ComputingTimeCase2}
\end{figure}

As depicted in Figure \ref{fig:ComputingTimeCase2}, the estimation time exhibits an interesting pattern. For instance, higher and more volatile computational times are observed for small values of $\lambda$, whereas a more regular and smoother pattern takes place when the value of $\lambda$ reaches the scale of generating units' costs (greater of equal to $10$ \$/MWh). The evaluation time exhibits a much more constant pattern as it does not directly depends on the value of $\lambda$.

\begin{figure}[h!]
    \centering
    \includegraphics[scale=.55]{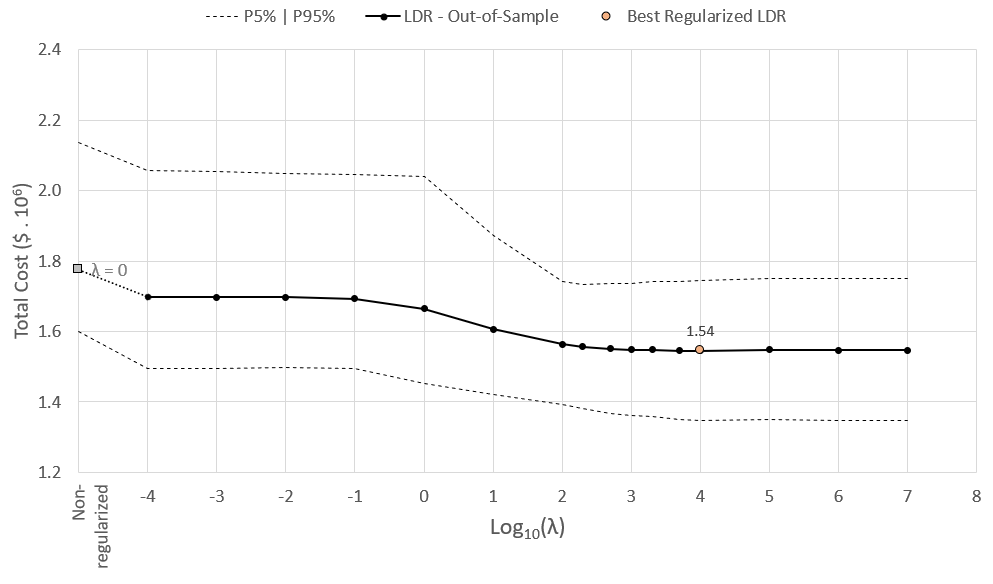}
    \caption{Expected Value of Total Operation Cost for different values of $\lambda$}
    \label{fig:Case2_TotalCost_Average}
\end{figure}

In this case study, we present the out-of-sample results for the expected value, fifth percentile (P5), and ninety-fifth percentile (P95) of the total operation cost. Figure \ref{fig:Case2_TotalCost_Average} presents these three indices (vertical axis) for different values of $\lambda$ (horizontal axis).

\noindent The lower and upper dashed curves show the P5 and P95 values, respectively, while the continuous curve depicts the expected value. In this case, the best value of $\lambda$ is equal to {$10000$}. Note that the P95 {metric} and the average costs {do not} agree in which should be the best value of $\lambda$. However, this is not necessarily true for all percentiles.

The best regularization produces a reduction of {$13.0\%$} on the total cost in comparison to the non-regularized benchmark. After this penalty level, the performance of the out-of-sample tests reaches worse results. The selection of coefficients resulted in 60.2\% of non-zero coefficients in non-regularized LDR against {2.0\%} for {$\lambda = 10000$}. The shrinkage also resulted in a relevant reduction, reaching {97.6\%} of reduction of the $\ell_1$-norm of the coefficient vector. 

\begin{figure}[h!]
    \centering
    \includegraphics[scale=.55]{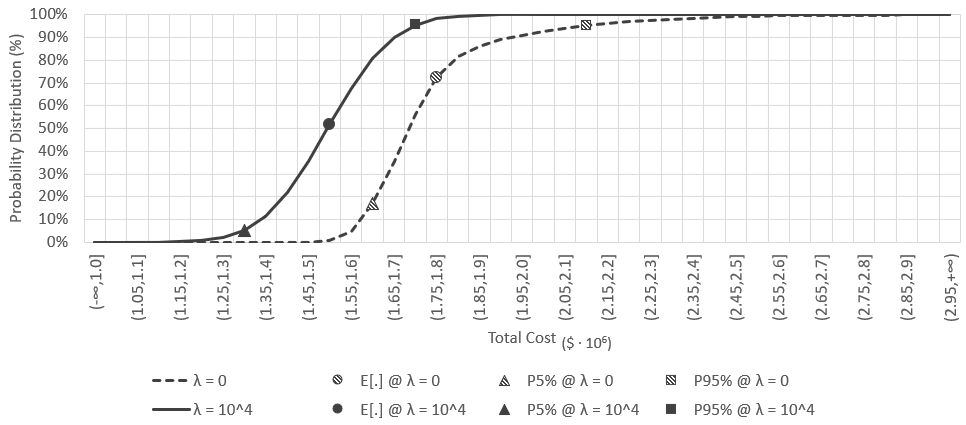}
    \caption{Distributed Probability Function, Expected Value, P5 and P95 of Operation Cost of Out-of-Sample Simulations}
    \label{fig:Case2_Out_Sample_ProbDist}
\end{figure}

In addition to the previously presented results, by comparing the accumulated probability distribution of both regularized and non-regularized LDRs, we can see in Figure \ref{fig:Case2_Out_Sample_ProbDist} that the regularized policy dominates (in the stochastic sense) the non-regularized one. In this figure, it is also highlighted the expected value, P5 and P95, also numerically presented in Tables \ref{tab:Case2_In_Sample_Results} and \ref{tab:Case2_Out-of-Sample_Results}. From the last column of Table \ref{tab:Case2_Out-of-Sample_Results}, we can also observe a 50\% reduction in the dispersion of the cost distribution (assessed through the difference between P95 and P5) for the regularized LDR case in comparison to the benchmark. Thus, the regularized LDR provides a stochastic dominant and less uncertain operational cost distribution. The previously highlighted findings showcase the prominent superiority of the regularized LDR over the non-regularized benchmark. 

Interestingly, by contrasting values presented in Tables \ref{tab:Case2_In_Sample_Results} and \ref{tab:Case2_Out-of-Sample_Results} we can see that: the small increments (of {1.5\%}, {0.6\%}, and {0.6\%}) on the in-sample metrics due to the regularization approach enables higher reductions (of {15.6\%}, {13.5\%}, and {18.7\%}) on the same indices in out-of-sample scenarios. This provides us with a measure of the in-sample sub-optimality imposed by the best regularization scheme.

\begin{table}[H]
    \begin{center}
    \caption{Cost metrics for in-sample scenarios.}
    \label{tab:Case2_In_Sample_Results}
    \begin{tabular}{c|c|c|c}
    \hline
    {$\lambda$} & {P5} & {Expected value} & {P95} \\
          & ($\$\cdot10^{6}$) & ($\$\cdot10^{6}$) & ($\$\cdot10^{6}$)\\
    \hline 
          0            & 1.34 & 1.58 & 1.77 \\
          {$10^{4}$} & {1.36} & 1.59 & 1.78 \\
    \hline
          Difference & {1.5\%} & {0.6\%} & {0.6\%} \\
    \hline
    \end{tabular}
    \end{center}
\end{table}
\begin{table}[H]
\begin{center}
\caption{Cost metrics for out-of-sample scenarios.}
\label{tab:Case2_Out-of-Sample_Results}
\begin{tabular}{c|c|c|c|c}
\hline
{$\lambda$} & {P5} & {Expected value} & {P95} & {P95}$-${P5}\\
      & ($\$\cdot10^{6}$) & ($\$\cdot10^{6}$) & ($\$\cdot10^{6}$) & ($\$\cdot10^{6}$)\\
\hline 
      $0$            & {1.60} & {1.78} & {2.14} & {0.54}\\
      {$10^{4}$} & {1.35} & {1.54} & {1.74} & {0.39}\\
    \hline
    Difference & {-15.6\%} & {-13.5\%} & {-18.7\%} & {-27.7}\%\\
    \hline
\end{tabular}
\end{center}
\end{table}

In view of evaluating the economic effects of the regularization, Table \ref{tab:Case2_MargCost} presents some additional important marked-oriented metrics that exploit the economic impact of the regularization process in spot prices (calculated in this work as the marginal operating costs). 
\begin{table}[H]
\begin{center}
\caption{Annual-average spot price metrics for out-of-sample scenarios.}
\label{tab:Case2_MargCost}
\begin{tabular}{c|c|c|c|c|c}
\hline
{$\lambda$} & {P5} & {E[.]} & {P95} & {Average} & {Time} \\
        & & & & {Uncertainty} & {Variability}
    \\
      & (\$/MWh) & (\$/MWh) & (\$/MWh) & (\$/MWh) & (\%)\\
\hline 
      $0$            & {340} & {622} & {1058} & {719} & {199}\\
      {$10^{4}$} & 148 & {195} & {314} & {166} & {27}\\
\hline
    \end{tabular}
\end{center}
\end{table}
\vspace{-0.5cm} 
Besides the usual metrics as percentiles and expected value for the average-annual spot prices (all extracted from the twelve central months of operation -- steady state), the average uncertainty level, calculated as the difference between the monthly P95 and P5, and the time variability, defined according to expression \eqref{TimeVar}, are also presented. 
\begin{equation} \label{TimeVar}
    \text{{Time Variability}} = \dfrac{1}{T-1}\sum_{\omega}{p_{\omega}\sum_{t=2}^{T}{\Big|\frac{\pi_{t, \omega}-\pi_{t-1, \omega}}{\pi_{t-1, \omega}}\Big|}}.
\end{equation}

Noticeably, the regularization policy also promotes remarkable stability in spot prices. Such an effect is characterized by a significant drop in the average uncertainty level and time variability indices when compared to the non-regularized policy. Additionally, the average of the annual spot price and its P95 also exhibit relevant reductions, which are desirable characteristics generally associated with efficiency signals for market players.

\subsubsection{Sensitivity Analysis}

{Finally, we present in Figure \ref{fig:Case2_RegBenefit&Time} a sensitivity analysis on the out-of-sample gain obtained with the regularized LDR (with the best value of $\lambda$) with respect to the number of scenarios $N$ used in the estimation process. It is possible to see that, although it is well known that the gain should decay to zero as $N$ (and $M$) grows to infinity, in practice, the size of tractable instances does not comport large sample sizes in the estimation step (see \cite{Bodur_2018} and other applications). In our case study, for $N \ge 500$, the model could not be solved and the computer run out of memory.

\begin{figure}[h]
    \centering
    \includegraphics[scale=.55]{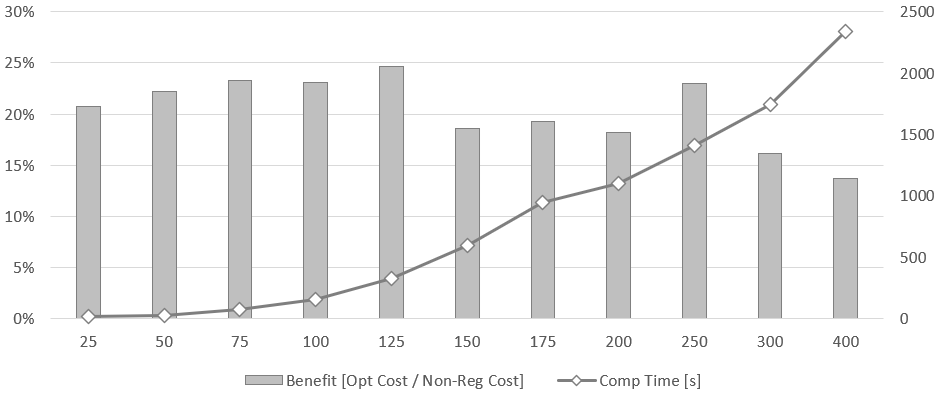}
    \caption{Sensitivity analysis of the regularization benefit in the out-of-sample with the estimation sample size $N$. Left axis (bars) -- out-of-sample gains obtained with the regularization procedure. Right axis (line) -- computational time took for the estimation process.}
    \label{fig:Case2_RegBenefit&Time}
\end{figure}

Although more powerful computers or decomposition algorithms could be implemented, in practical applications, SAA instances are in general constrained to finite (small) sample sizes. This is especially true in real industry applications, where we need to consider very detailed models to avoid unrealistic (and even infeasible) solutions.\footnote{For instance, the official dispatch planning models used in Brazil use very few scenarios (less than 40). Planners also rely on models with very few scenarios due to the huge number of constraints to model physical aspects needed to be considered to ensure feasible solutions (see \cite{Soares2022} for a realistic planning study based on a new decomposition algorithm and detailed comments on the subject).} Notwithstanding, previously reported work on two-stage LDR also rely on small sample sizes (see \cite{Bodur_2018} and \cite{daryalal2020lagrangian}). In this context, we can see in Figure \ref{fig:Case2_RegBenefit&Time} that for all tractable instances analyzed in this case study, the proposed regularized LDR improved the SAA result by more than 10\%\footnote{{All reported gains are statistically significant at a significance level lower than 0.01.}}, which is a considerable amount for practical applications. \emph{Therefore, the above results, although empirical and problem based, provide large evidence in favor of corroborating the initial hypothesis, namely, that LDR-based policies estimated with AdaLASSO-regularization can provide relevant cost savings under out-of-sample (unseen) data.}}



\section{Conclusion} \label{Conclusion}

This work raised awareness of the in-sample overfit threat featured by linear decision rules (LDRs) applied to solve multistage stochastic linear problems (MSLP). In this context, we proposed a new regularized LDR-based policy considering the AdaLASSO, which takes into account the tradeoff between a better in-sample fit and the over parameterization of the LDR. Thus, we tested the following hypothesis: \emph{non-regularized LDR overfits in-sample data and more parsimonious LDR-based policies estimated with AdaLASSO-regularization can provide relevant cost savings under unseen out-of-sample data.} We found strong evidence that this hypothesis is true based on our computational results. 

To study the proposed framework and test the aforementioned hypothesis, we used the well-known and representative long-term hydrothermal dispatch planning problem (see \cite{Pereira1991}) considering a basis of nonlinear functions in the two-stage LDR (see \cite{Bodur_2018}). {Within the limitations of our case study (selected problem and data),} results show that, in the out-of-sample test, the regularized LDR achieved relevant cost savings in comparison to the classical non-regularized LDR approach (18.5\% on average and 25.0\% for the 95-percentile in a single-reservoir case study and {13.5\%} on average and {18.7\%} for the 95-percentile of the costs in a multi-reservoir case study). The AdaLASSO regularization resulted in a reduction greater or equal to 82.9\% of the total number of selected coefficients in both case studies. This fact highlights the importance of the regularization scheme to find the subset of parameters that performs best in out-of-sample. Interestingly, our case study also revealed a relevant effect of the best regularization on dual variables. The uncertainty level and time variability metrics of spot prices (marginal operating costs extracted from a steady state year) are significantly reduced when compared to the non-regularized policy. The expected value and the 95 percentile of the annual--average spot prices also show a relevant reduction under the best regularization penalty. These facts provide strong evidence that the proposed regularization approach constitutes an important step, worth analyzing when using LDRs to solve MSLP. 

{ Notwithstanding, relevant and interesting challenges arise from the utilization of LASSO-based regularization terms to address the overfitting issue when using LDRs to address MSLP. For instance, the identification of the best penalty parameter, $\lambda$, requires a line search procedure running the MSLP for each point inspected. Interestingly, results suggest that the best parameter is stable across instances, which is a merit of this approach as one might not need to calibrate the parameter every time. Additionally, the use of sensitivity analysis can help to identify regions for which the LDR remains unchanged and thereby a recalculation is not necessary.}

We highlight {four} possible avenues for future research in this topic: 1) the study of theoretic results based on convergence results for two-stage sample average approximation such as \citep{ahmed2002sample} and \citep{kleywegt2002sample}; 2) the study of different applications in which traditional methods' hypothesis, such as convexity in the case of the SDDP, fail to comply with the problem at hand, e.g., non-convex (nonlinear) hydro production curves, unit-commitment (binary) constraints, revenue maximization with spot-price uncertainty, just to mention a few; 3) the study of different bases of functions and regularization methods, e.g., $\ell_2$-norm {and combinations between $\ell_1$ and $\ell_2$-norm (which are associated with ridge regression and Elastic Net) with interesting and complementary properties to AdaLASSO}; finally, 4) {the study of interactions between scenario reduction techniques and regularization.}

\bibliography{mybibfile}

\end{document}